%% file: main.tex
\documentclass[11pt]{article}
	\newcommand{\blind}{0}
	
    \makeatletter
    \renewcommand\section{\@startsection {section}{1}{\z@}%
                                       {-3.5ex \@plus -1ex \@minus -.2ex}%
                                       {2.3ex \@plus.2ex}%
                                       {\normalfont\fontfamily{phv}\fontsize{16}{19}\bfseries}}
    \renewcommand\subsection{\@startsection{subsection}{2}{\z@}%
                                         {-3.25ex\@plus -1ex \@minus -.2ex}%
                                         {1.5ex \@plus .2ex}%
                                         {\normalfont\fontfamily{phv}\fontsize{14}{17}\bfseries}}
    \renewcommand\subsubsection{\@startsection{subsubsection}{3}{\z@}%
                                        {-3.25ex\@plus -1ex \@minus -.2ex}%
                                         {1.5ex \@plus .2ex}%
                                         {\normalfont\normalsize\fontfamily{phv}\fontsize{14}{17}\selectfont}}
    \makeatother
	
	\usepackage[utf8]{inputenc}
    \usepackage[T1]{fontenc}
    \usepackage{amsmath,amssymb}
    \usepackage{lipsum}
    \usepackage[short,nocomma]{optidef}
    \usepackage[margin=1.00in]{geometry}
    \usepackage{enumitem}
    \usepackage{multirow,multicol}
    \usepackage[hidelinks]{hyperref}
    \usepackage{algorithm}
    \usepackage{algpseudocode}
    \usepackage{xcolor}
    \usepackage{pdfpages}
    \usepackage{svg}
    \usepackage{caption}
    \usepackage{scalerel}

\DeclareMathOperator*{\minimize}{min}
\DeclareMathOperator*{\maximize}{max}

\def\tigerdams{Tiger Dams\texttrademark~}

\begin{document}

\def\spacingset#1{\renewcommand{\baselinestretch}%
    {#1}\small\normalsize} \spacingset{1}

\if0\blind
{
    \title{\bf Scenario-based Optimization Models for Power Grid Resilience to Extreme Flooding Events}
    \author{Ashutosh Shukla, Erhan Kutanoglu, John J. Hasenbein\\
    Graduate Program in Operations Research and Industrial Engineering\\ The University of Texas at Austin, Austin, United States}
    \date{}
    \maketitle
} \fi

\if1\blind
{

    \title{\bf Scenario-based Optimization Models for Power Grid Resilience to Extreme Flooding Events}
    \author{Author information is purposely removed for double-blind review}
    \maketitle
 }  \fi 
		
\begin{abstract}
We propose two scenario-based optimization models for power grid resilience decision making that integrate output from a hydrology model with a power flow model. The models are used to identify an optimal substation hardening strategy against potential flooding from storms for a given investment budget, which if implemented enhances the resilience of the power grid, minimizing the power demand that is shed. The same models can alternatively be used to determine the optimal budget that should be allocated for substation hardening when long-term forecasts of storm frequency and impact (specifically restoration times) are available. The two optimization models differ in terms of capturing risk attitude: one minimizes the average load shed for given scenario probabilities and the other minimizes the worst-case load shed without needing  scenario probabilities. To demonstrate the efficacy of the proposed models, we further develop a case study for the Texas Gulf Coast using storm surge maps developed by the National Oceanic and Atmospheric Administration and a synthetic power grid for the state of Texas developed as part of an ARPA-E project. For a reasonable choice of parameters, we show that a scenario-based representation of uncertainty can offer a significant improvement in minimizing load shed as compared to using point estimates or average flood values. We further show that when the available investment budget is relatively high,  solutions that minimize the worst-case load shed can offer several advantages as compared to solutions obtained from minimizing the average load shed. Lastly, we show that even for relatively low values of load loss and short post-hurricane power restoration times, it is optimal to make significant investments in substation hardening to deal with the storm surge considered in the NOAA flood scenarios.

\end{abstract}
			
\noindent%

\spacingset{1.5} 

\input{introduction}
\input{literature}
\input{model}
\input{case}
\input{results.tex}


\bibliographystyle{unsrt}
\spacingset{1}
\bibliography{main}


\end{document}

%% file: introduction.tex
\section{Introduction}

In the past few years, hurricanes and tropical storms have caused significant damage to critical infrastructures such as transportation systems, healthcare services, and the power grid. Hurricanes Maria, Irma, and Harvey together cost nearly \$265B which was more than 85\% of total weather-related disaster costs in the U.S. in 2017 \cite{billion_nerc}. Harvey became not only the longest-lasting hurricane with a record level of rainfall but also the costliest at \$130B, part of which was due to power outages. Harvey damaged 90+ substations, downed 800+ transmission assets, 6000+ distribution poles, and 800+ miles of power lines, with a peak power generation loss of 11GW, affecting over 2 million people. It took 2 weeks and 12,000 crew members to restore power \cite{NERC2018}. 

The power grid is impacted by hurricanes and tropical storms primarily due to strong winds and flooding. To address this, a vast body of research that examines the effect of wind fields on transmission lines and towers has been developed. However, to the best of our knowledge, the literature is quite scant on the models that assess the impact of flooding. At the same time, the cost of such disasters has increased in states like Texas which is exposed to the Atlantic basin through the Gulf of Mexico. Moreover, recent studies suggest that we are likely to see more frequent and intense hurricanes in the near future \cite{Webster2005}. In response to this, some utilities have employed on-site meteorologists which they have reported to be beneficial \cite{NERC2018}. These meteorologists localize the predictions to obtain flood estimates for the region of interest. The estimates are then used to determine the resources needed for forecasted damage and post-storm recovery. To further improve this  decision-making process, we present an end-to-end scenario-based optimization approach that integrates the output from a predictive geoscience-based flood model with a power flow model to recommend a plan for substation hardening to relieve the flood impacts of the potential storms. While doing so, the scenario-based approach accounts for the uncertainty associated with storms and their flood forecasts. 

Specifically, we propose two scenario-based optimization models (stochastic and robust) for grid resilience decision making under uncertainty. The choice of the model to be used for decision-making depends on the available information about the uncertain parameters which in our case are the flood levels at substations in a flood scenario. We show how the proposed models can be used to identify the substations that should be protected and to what extent. We further explain how the same models can be used for deciding the optimal budget that should be allocated for substation hardening to minimize the expected total disaster management cost. The aforementioned features of the models can help power utilities and grid operators address their concerns like the unpredictable nature of the load loss, the potential for substation flooding, and the potential reduction in generator output due to loss of load as outlined in \cite{NERC2018}.

The rest of the article is organized as follows. Section \ref{lit_review} presents a review of the literature on power grid resilience, particularly from a modeling and decision making viewpoint. Section \ref{model} presents the overview of the proposed models followed by the notation, assumptions, mathematical formulation, and a brief discussion on the characteristics of the models. Section \ref{case_study} is dedicated to the development of a case study for the Texas Gulf Coast and Section \ref{results} is to the discussion of the results. We conclude with directions for future research in Section \ref{conclusions}.

%% file: literature.tex
\section{Literature Review} \label{lit_review}

Power grid resilience to extreme events like cyber-attacks and natural disasters has been a topic of intense research in the past few years \cite{framework_grid_resilience}. This includes studies focused on developing resilience metrics, methodological frameworks to enhance power grid resilience, and approaches to risk analysis \cite{Panteli_framework, Panteli_metrics, framework_5}. In addition, several mathematical models have been developed to aid decision-making in different stages of the power grid resilience management cycle. These models can be categorized based on the planning phase they are developed for: mitigation, preparedness, response, and recovery. Using the flood as an extreme event that the resilience models are designed to respond to, the mitigation decisions are about the permanent hardening of the grid components, re-design of the grid through the introduction of new substations/transmission lines, installation of backup generation, etc. These decisions are made well before the start of the hurricane season when limited information about upcoming hurricanes is available. Similarly, before an imminent hurricane, during the preparedness phase, one has to make decisions about where to install temporary flood barriers like \tigerdams, where to deploy mobile substations to quickly recover damaged substations, and what part of the grid to disconnect to avoid fatal accidents due to the collapse of power lines. In both the mitigation and preparation phases, decision-makers face significant uncertainty about storm characteristics like path, intensity, forward speed, and precipitation.
The decision-making in the response and recovery phase, on the other hand, is not plagued by weather uncertainty.
Since this paper focuses on decision-making under weather uncertainty, we limit the review's focus to models that aid decision-making during the mitigation and preparedness phases. 

Models for both the mitigation and preparedness phases can be broadly categorized into two groups: (1) machine learning-based and (2) optimization-based. In the case of machine learning-based models, the focus is on prediction, not decision making. For example, a machine learning-based model may predict metrics of interest such as the number of outages, outage duration, etc., for an upcoming hurricane \cite{seth_ml_outage_prediction, seth_ml_outage_prediction_2, Reilly2015}. However, decision making based on these predictions, like which substations to protect and how to reconfigure the power grid network to minimize load shed, are not typically considered within the model. 
Optimization-based models leverage predictions for decision making. To do so, the predictions with associated uncertainty are represented using scenarios. The decision-making model is then coupled with these scenarios. The models that we propose in this study belong to the class of optimization-based models. In the subsequent paragraphs, we survey the key characteristics of some of these models and highlight their differences from what we propose.

The optimization models generally consist of two components: uncertainty quantification and decision modeling. To quantify uncertainty about the weather, we first generate a set of scenarios using various kinds of models, such as machine learning-based models, physics-based models, and expert opinions. Then, irrespective of how the scenarios are generated, the decision-making model considers the impacts on the grid under each scenario to recommend decisions that minimize a certain risk measure. The models we review in this subsection are based on different ways the aforementioned components of the optimization model can be developed. In particular, we review various methods of generating representative scenarios and incorporating them into alternative decision making models. 

\subsection{Scenario generation}

Scenario generation is one of the most common uncertainty quantification methods for extreme weather. We divide scenario generation techniques into four categories. The first is based on fragility curves. The curve represents the failure probability of a component as a function of some loading parameter. For example, fragility curves for transmission towers have been developed with respect to wind speeds. Such curves have been used in various power grid resilience decision making studies \cite{seth_2021, Panteli_boosting}. The second is based on statistical methods. For example, in \cite{Staid2014}, the authors use historical hurricane and tropical storm data for developing a baseline scenario. The alternative scenarios are then developed by altering parameters from the historical data to simulate plausible climate-induced changes to storm behavior. 
In \cite{9842364}, the path and the wind field of typhoons are simulated using Monte Carlo sampling to quantify the spatio-temporal impacts of wind speed on the transmission line status. In addition to wind and flooding, winter storms in Texas, such as Uri of 2021, have propelled research in power grid resilience to extreme cold events. For example, in \cite{austgen_winter_2022}, the authors have developed a statistical model where they incorporate historical outage data to generate scenarios of generator outages due to extreme cold events. The third set of methods is based on physics-based hydrological models. Two such models, called WRF-Hydro and SLOSH, are used in \cite{kyoung_paper, SOUTO2022107545} to generate flooding scenarios. In \cite{9695371}, the authors use physics-based climate models to evaluate the resilience of levee-protected electric power networks with the primary focus on performance degradation. The fourth category is based on combinatorial criteria like $N-k$. In this case, each scenario represents a way in which $k$ out of $N$ components can fail. A model based on this criterion is used in \cite{9470975}.

\subsection{Decision modeling} \label{decision_model}

Several optimization models have been developed for power grid resilience decision-making against extreme weather events. These include models that can aid in decision-making about the upgrade of the power grid network through a combination of hardening existing components, adding redundant lines, switches, generators, and transformers \cite{9842364, 9470975, 7540988}. However, hardening large parts of the power grid can be financially infeasible. In such a scenario, stockpiling power grid components in strategic locations enhances resilience by expediting network restoration after the disaster. To decide how stockpiling of components should be done, Coffrin \textit{et al}. \cite{Coffrin2011} developed a two-stage stochastic mixed-integer program where the first-stage discrete decisions are about stockpiling power grid components and the second-stage decisions are about how to operate the power grid to minimize load-shed. 
Additionally, network reconfiguration before an imminent hurricane can also enhance resilience. The models proposed in \cite{Panteli_boosting} make such decisions using grid islanding techniques. However, none of these models have explicitly focused on assessing the impact of flooding on the power grid. On the other hand, there are several studies that assess the impact of flooding on other critical infrastructures.
For example, Kim et al. \cite{kyoung_paper} present a framework and a case study using hurricane Harvey to generate physics-based (hydrological) flood scenarios. These scenarios are then used for resilience decision-making for healthcare infrastructure in \cite{gizem_2022}. Scenarios generated from physics-based models have also been used in \cite{Yuki_Columbia_model} that developed a model to estimate the overall disaster cost due to physical damage loss, income losses, and inventory losses. In comparison, our proposed models are explicitly geared towards resilience decision-making for the power grid and have a power flow model nested within a larger substation hardening model. 

The models closest to ours are \cite{SOUTO2022107545} and \cite{Mohadese}. \cite{SOUTO2022107545} uses a set of scenarios all based on Hurricane Harvey run on a hydrology model focusing on inland flooding. We instead consider a wider range of storms and storm characteristics and scenarios that are based on NOAA's storm surge simulations. Mohadese et al. \cite{Mohadese} propose a stochastic optimization model for identifying and protecting substations a day before the anticipated flooding event, meaning a focus on preparedness. 
Here, our proposed models differ in several ways. First, we generate scenarios using outputs from physics-based hydrological models to create flood maps for the region of interest. Our choice is based on the rationale that physics-based models represent flood levels across the region of interest such that they are correlated in space and time. Mohadese et al. \cite{Mohadese} have not considered the impact of correlated flooding. They also assume that a substation will transmit power if it is not flooded. In reality, this may not be true due to the network effects and we embed a power flow model within our larger resilience optimization model to address such effects. Finally, our models focus on long-term decision making, highlighting mitigation-phase budgeting and decision making. 

%% file: model.tex
\section{Modeling} \label{model}

In this section, we first present an overview of the stochastic and robust optimization models developed to assist in grid resilience decision-making against extreme flooding events. Next, we state the key assumptions, introduce the notation and provide detailed mathematical formulations. Finally, we highlight some of the key characteristics of the proposed models and explain how they can be used to address a wide variety of questions in grid resilience decision-making. We highlight that the models proposed in subsection \ref{SO_model} and \ref{RO_model} are developed to minimize a risk measure over the load shed for a single flood event. In subsection \ref{model_discussion}, we explain how the same models can be leveraged for multi-year planning.

\subsection{Overview of the proposed models}\label{overview}

\begin{figure}
\captionsetup{width=0.45\textwidth}
\centering
\begin{minipage}{.45\textwidth}
    \centering 
    \includegraphics[width=0.95\textwidth]
    {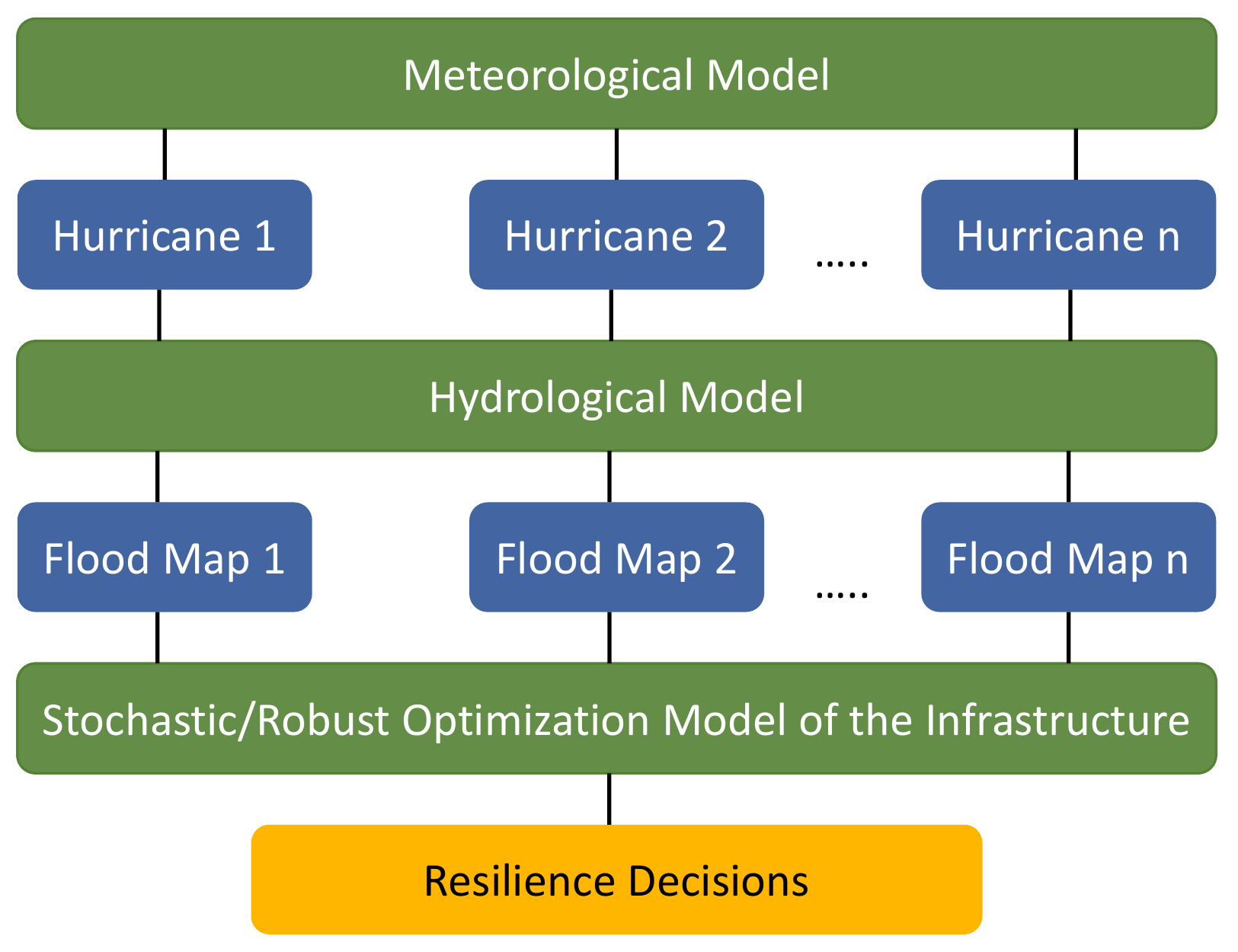}
    \captionof{figure}{A schematic representation of the decision-making framework for the proposed models}
    \label{fig:schematic_only}
\end{minipage}%
\begin{minipage}{.55\textwidth}
    \centering
    \vstretch{1.25}{\includegraphics[width=3.3in]{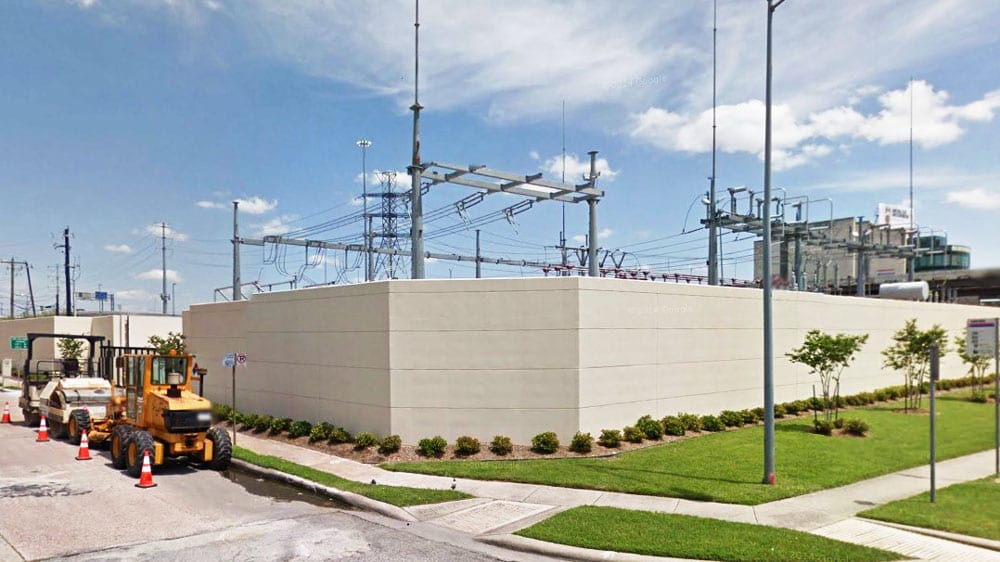}}
    \captionof{figure}{An example permanent hardening structure at a substation}
    \label{fig: permanent}
\end{minipage}
\end{figure}

A two-stage stochastic optimization model is developed to address situations where the uncertainty about hurricane-induced flooding is modeled using a probability distribution. In this case, the model minimizes the expected unsatisfied power demand (also referred to as the expected load shed) due to the components' failures (i.e., flooded substations) over a set of scenarios. The two-stage robust model on the other hand requires no information about the probability distribution. The model instead minimizes the maximum load shed in any scenario within the uncertainty set. A general framework representative of the proposed models is shown in Figure \ref{fig:schematic_only}. As shown in the figure, we represent the uncertainty in the meteorological forecasts using a set of hurricane parameters (Hurricane $1, \ldots, n$). In the next step, using the aforementioned parameters as input, we run a hydrological model to get the corresponding flood maps. The flood maps are then used as input to the two-stage decision-making models. The final output from the decision-making models is a plan for substation hardening.

A scenario in the proposed models represents the water levels at different substations obtained from the flood map of a specific hurricane type. Hurricanes with different characteristics such as direction, intensity, forward speed, etc.\ lead to different levels of flooding, generating different scenarios. These scenarios are representative of the flooding that the region of study can experience over a specific time period (typically multiple years). A distinguishing feature of the proposed models is the way we generate the scenarios. Instead of using popular techniques like fragility curves, we use flood maps for scenario representation because they capture the effects of correlated flooding. This is important because the failure of a substation within a power grid can have network effects on the other parts of the grid. To account for such details, we need not only know which substations fail more frequently but also the combination of substations that fail together. Our proposed model accounts for these details and uses them to evaluate the effects of such damages (in the form of load shed) during decision-making by solving a power flow model. 

The power grid network considered in the proposed models is represented by a graph where the buses are represented by the nodes and the branches interconnecting the buses are represented by the edges of the graph. The branches of the network are held using the transmission towers. We assume that these towers are well above the ground
and therefore immune to flooding. It is the substations and components within them that are susceptible to flooding. In this study, we assume that the substations are outdoor open-air facilities. Therefore, when a substation is flooded, we assume that all the components within the substation and the branches connected to all the buses within the substation are out of order. 


Decision making in both models occurs in two stages. Here we specify first-stage decisions that determine which substations to harden and to what extent. We assume that the substation hardening measures are taken during the mitigation phase  of the power grid resilience management cycle. Consequently, the decisions made using the model are not a response to any particular imminent hurricane. Instead, they are intended to harden the grid against multiple hurricanes potentially occurring over multiple years and minimize the long-term disaster costs incurred due to flooding. One such mitigation measure for substation hardening is to build permanent protective structures like walls around the substation periphery as shown in Figure \ref{fig: permanent}. 

After we make the first-stage decisions, we assess their performance in dealing with the flood levels in the second stage. The second-stage assessment involves minimization of the load shed in multiple flood scenarios that may impact the power grid during the multi-year planning horizon. For each of these scenarios, we overlay the power grid network on the flood map to identify parts of the network that are flooded. Given the flood height and the level of hardening at a substation, the model infers if a substation is flooded in a particular scenario. If a substation is flooded, all the buses within the substation and the branches connected to those buses are considered to be out of order. Once the damaged state of the power grid network is determined, we solve the second-stage assessment (the so-called recourse problem) which is a DC power flow model to estimate the load shed given the state of the grid. The second-stage decisions in the recourse problem determine the routing of power to minimize unsatisfied power demand. It should be noted that although both the stochastic and robust models involve the same sets of decision variables, the specific solutions suggested by them can be vastly different. The robust model gives us the flexibility to make decisions in the absence of any information about the probability distribution. These decisions can however be far more conservative than the decisions recommended by the stochastic model.


\subsection{Assumptions}

In the proposed models, we make the following assumptions. The first assumption is that a DC-based power flow approximation is acceptable. For a detailed explanation and derivation of the DC power flow equations from the AC equations, we refer to \cite{8635446}. This approximation has been widely used and is embedded within larger strategic decision-making problems such as long-term capacity planning and operation of wholesale electricity grids. For the kind of models proposed in this paper, a detailed discussion on the difference in the quality of solutions obtained from different power flow approximation models is given in \cite{austgen_impacts_2021}. The second assumption is that we can model the substation hardening cost with fixed and variable components. The fixed cost is incurred when a substation is chosen for hardening. It can represent the cost of building the foundation on which the protective structure is built. Furthermore, it can also include the costs associated with transporting construction resources to the substation site. The variable cost, on the other hand, is a function of the height of flooding to which the substation is made resilient. In our case, we assume that the variable cost linearly depends on height. This is reasonable when we build wall-like structures to protect substations, as shown in Figure \ref{fig: permanent}. Third, we assume that each substation's hardening and flood levels are discrete and finite. In the proposed formulations, they are assumed to be non-negative integers. Fourth, we assume that all the flooded substations within the network experience the same downtime and are recovered simultaneously. Lastly, we assume that the value of load loss can be quantified in dollars per hour.

\subsection{Notation}

The proposed models use the following notation. Note that all the cost parameters used in the models are in dollars and all power grid parameters are in the per-unit system. Notation not defined in this section and appearing later in the text is defined as introduced.

\noindent\textbf{Sets}
\begin{itemize}[noitemsep,nolistsep]
    \item[] $\mathcal{I}$: Set of substations indexed by $i$
    \item[] $\mathcal{I}_f$: Set of substations that are flooded in at least one scenario
    \item[] $\mathcal{J}$: Set of buses indexed by $j$
    \item[] $\mathcal{K}$: Set of scenarios indexed by $k$
    \item[] $\mathcal{R}$: Set of branches indexed by $r$
    \item[] $\mathcal{B}_i$: Set of buses at substation $i$
    \item[] $\mathcal{N}_j^{in}$: Set of branches incident on bus $j$ with power flowing into bus $j$
    \item[] $\mathcal{N}_j^{out}$: Set of branches incident on bus $j$ with power flowing out of bus $j$
\end{itemize}

\noindent\textbf{Parameters}
\begin{itemize}[noitemsep,nolistsep]
    \item[] $M$: An arbitrarily large constant
    \item[] $f_i$: Fixed cost of hardening at substation $i$
    \item[] $v_i$: Variable cost of hardening at substation $i$
    \item[] $H_i$: Maximum flood height to which substation $i$ can be hardened
    \item[] $\theta(j)$: Substation that contains bus $j$
    \item[] $\Delta_{ik}$: Flood height at substation $i$ in scenario $k$ (a non-negative integer value)
    
    \item[] $B_r$: Susceptance of branch $r$
    \item[] $F_r$: Maximum power that can flow in branch $r$
    
    \item[] $\lambda(r)$, $\mu(r)$: Head bus and tail bus of branch $r$
    \item[] $D_j$: Load at bus $j$
    \item[] $\underline{G}_j$, $\overline{G}_j$: Minimum and maximum generation at bus $j$
    \item[] $\beta$: Index of the reference bus
    \item[] $p_k$: Probability of scenario $k$
    \item[] $I$: Total investment budget for substation hardening
\end{itemize}

\noindent\textbf{Variables}
\begin{itemize}[noitemsep,nolistsep]
    \item[] $y_i$: Binary variable indicating whether substation $i$ is chosen for permanent hardening
    \item[] $x_i$: Non-negative integer variable indicating discrete height of hardening at substation $i$
    \item[] $z_{j}$: Binary variable indicating if bus $j$ is operational
    \item[] $s_{j}$: Non-negative real variable indicating load satisfied at bus $j$
    \item[] $g_{j}$: Non-negative real variable indicating power generated at bus $j$
    \item[] $u_{j}$: Binary variable indicating if generator at bus $j$ is used
    \item[] $\alpha_{j}$: Real variable indicating voltage phase angle of bus $j$
    \item[] $e_{r}$: Real variable indicating power flowing in branch $r$
\end{itemize}

\subsection{Stochastic Optimization Model}\label{SO_model}

The two-stage stochastic optimization model (SO) is expressed as:
\begin{subequations}
\begin{equation}\label{SO_problem}
    L_{\mathcal{SO}}^* = \displaystyle{\minimize_{x \in \mathcal{X}}} \, \, L_{\mathcal{SO}}(x),  
\end{equation}
\text{where} 
\begin{equation}\label{SO_defination}
    L_{\mathcal{SO}}(x) = \sum_{k \in \mathcal{K}} p_k\mathcal{L}(x,k).
\end{equation}
\end{subequations}
The objective function in \eqref{SO_problem} minimizes the expected unsatisfied power demand (load shed) over the scenarios in set $\mathcal{K}$. 
Here, $\mathcal{X}$ represents the set of feasible first-stage decisions. The following constraints define the set:
\begin{subequations}
\begin{equation}\label{budget_constraint}
        \sum_{i \in \mathcal{I}_f} f_i y_{i} + v_i x_{i} \leq I,
    \end{equation}
    \begin{equation}\label{constraint_tighetening}
       x_{i} \leq H_i y_i, \quad \forall  i \in \mathcal{I}_f. 
    \end{equation}
\end{subequations}
Note that the variables $x_i$ and $y_i$ are defined only for substations that are flooded in at least one scenario. Constraint (\ref{budget_constraint}) enforces that the sum of the fixed and variable costs incurred due to substation hardening does not exceed the investment budget. Constraints (\ref{constraint_tighetening}) place an upper bound on the extent of flooding to which the substation can be made resilient while linking variables $x_i$ and $y_i$ for each substation $i$. Such constraints represent engineering and practical challenges that may arise while building protective structures that are too tall.

In \eqref{SO_defination}, $L_{\mathcal{SO}}(x)$ represents the expected load shed when the first-stage decision is $x$. Here, $p_k$ represents the probability of scenario $k$ and $\mathcal{L}(x,k)$ is the recourse function representing the minimum load shed when the first-stage decision is $x$ and scenario $k$ is realized. The recourse function is defined as follows:
\begin{mini!}|l|[1]<b>
    {}{\sum_{j \in \mathcal{J}} D_j - s_{j},
       \label{objective}}{}{\mathcal{L}(x,k) =}{}
    \addConstraint{(1-z_{j})M \geq \Delta_{\theta(j)k} - x_{\theta(j)},}
                  {}{\forall j: \theta(j) \in \mathcal{I}_f,}
                  {\label{linking_1}}
    \addConstraint{2 z_{j}M \geq 1 - 2(\Delta_{\theta(j)k} - x_{\theta(j)}),}
                  {}{\forall j: \theta(j) \in \mathcal{I}_f,}
                  {\label{linking_2}}
    \addConstraint{z_j = 1,}
                  {}{\forall j: \theta(j) \in \mathcal{I} \setminus \mathcal{I}_f,}
                  {\label{other_z}}
    \addConstraint{u_{j} \leq z_{j},}
                  {}{\forall j \in \mathcal{J},} {\label{flexible_generation}}
    \addConstraint{s_{j} \leq z_{j}D_j,}
                  {}{\forall j \in \mathcal{J},}
                  {\label{supply_demand}}
    \addConstraint{u_{j}\underline{G}_j \leq g_{j} \leq u_{j}\overline{G}_j,}
                  {}{\forall j \in \mathcal{J},}
                  {\label{capacity_constraint}}
    \addConstraint{-z_{\lambda(r)} F_r \leq e_{r} \leq z_{\lambda(r)} F_r,}
                  {}{\forall r \in \mathcal{R},}
                  {\label{edge_operational_1}}
    \addConstraint{-z_{\mu(r)} F_r \leq e_{r} \leq z_{\mu(r)} F_r,}
                  {}{\forall r \in \mathcal{R},}
                  {\label{edge_operational_2}}
    \addConstraint{B_r (\alpha_{\lambda(r)} - \alpha_{\mu(r)})
                   \geq M(z_{\lambda(r)} + z_{\mu(r)}) - 2M + e_{r},}
                  {}{\forall r \in \mathcal{R},}
                  {\label{phase_angle_constraint1}}
    \addConstraint{B_r (\alpha_{\lambda(r)} - \alpha_{\mu(r)})
                   \leq -M(z_{\lambda(r)} + z_{\mu(r)}) + 2M + e_{r}, \quad}
                  {}{\forall r \in \mathcal{R},}
                  {\label{phase_angle_constraint2}}
    \addConstraint{\sum_{r \in \mathcal{N}_j^{out}}e_{r} - \sum_{r \in \mathcal{N}_j^{in}}e_{r}= g_{j} - s_{j},}
                  {}{\forall j \in J,}
                  {\label{flow_balance}}
    \addConstraint{-\pi \leq \alpha_{j} \leq \pi,}
                  {}{\forall j \in \mathcal{J},}
                  {\label{phase_side_constraint}}
    \addConstraint{\alpha_{\beta } = 0.}
                  {}{}
                  {\label{reference_bus}}
\end{mini!}

The objective function in (\ref{objective}) minimizes the unsatisfied power demand when the first-stage decision is $x$ and the flood scenario realized is $k$. Constraints (\ref{linking_1}) and (\ref{linking_2}) link the first-stage substation hardening decisions to the second-stage scenario-dependent power flow decisions. For a given hardening decision at a substation, the provided protection level is compared against the flood height
at that substation in a given scenario. Depending on whether the hardening level can withstand the flooding, we set the status of the corresponding bus as operational or not. This is indicated by variable $z_{j}$. For the substations that are not flooded in any of the scenarios, the status of the corresponding buses is set to operational in constraints \eqref{other_z}.

Constraints (\ref{flexible_generation}) capture generator dispatch decisions for operational generators. When $z_j = 0$, we cannot inject power to the network through bus $j$ and therefore set $u_{j} = 0$. If $z_j = 1$, we let the recourse problem decide if power generated at bus $j$ should be used or not. Constraints (\ref{supply_demand}) place an upper bound on the amount of power that can be supplied at bus $j$ (which is the demand at that bus). If bus $j$ is flooded, then $z_{j} = 0$ and no power can be supplied to the loads that are connected to the bus. Constraints (\ref{capacity_constraint}) place upper and lower bounds on the amount of power that can be generated at bus $j$. If bus $j$ is flooded, then $u_{j} = 0$ and thus $g_{j} = 0$. If on the other hand bus $j$ is not flooded, the model solves the recourse problem which is a binary linear program to determine the amount of power that should be generated at bus $j$. Constraints (\ref{edge_operational_1}) and (\ref{edge_operational_2}) place restrictions on the amount of power that can flow through branch $r$. If the bus at either end of a branch is flooded, then no power can flow through it. On the other hand, if buses at both ends of the branch are operational, then a maximum of $F_r$ power can flow through it in either direction. Constraints (\ref{phase_angle_constraint1}) and (\ref{phase_angle_constraint2}) enforce an approximation to Ohm's Law. If both ends of a branch are operational, then the amount of power flowing in the branch is governed by equations
$$B_r (\alpha_{\lambda(r)} - \alpha_{\mu(r)}) = e_{r}, \quad \forall r \in \mathcal{R}.$$
If the bus at either end of the branch is flooded, then the above equation need not hold. This is achieved by introducing big-$M$ values. The formulation can be further tightened by appropriately determining the values of big-$M$. A discussion on this is presented in Section \ref{model_discussion}. Constraints (\ref{flow_balance}) represent the flow balance which states that the net power injected into the network at bus $j$ is the difference between the power generated and consumed at the same bus. Constraints (\ref{phase_side_constraint}) impose limits on the phase angle values at buses. Finally, constraints (\ref{reference_bus}) set the phase angle of the reference/slack bus to 0. 

\subsection{Robust optimization model}\label{RO_model}

In the two-stage robust optimization model (RO), we minimize the maximum unsatisfied power demand value across all scenarios. Mathematically, the problem can be stated as follows:
\begin{subequations}
\begin{equation}\label{RO_problem}
    L_{\mathcal{RO}}^* = \displaystyle{\minimize_{x \in \mathcal{X}}} \, \, L_{\mathcal{RO}}(x),  
\end{equation}
\text{where} 
\begin{equation}\label{RO_defination}
    L_{\mathcal{RO}}(x) = \displaystyle{\maximize_{k \in \mathcal{K}}} \, \, \mathcal{L}(x,k).
\end{equation}
\end{subequations}
The expression in \eqref{RO_defination} finds the maximum scenario-based load shed via a max($\cdot$) function. RO in \eqref{RO_problem} can be reformulated as\begin{equation}\label{epigraph_constraint}
\displaystyle{\minimize_{x \in \mathcal{X}}} \, \, \left\{\tau : \tau \geq \mathcal{L}(x,k) \quad \forall \quad k \in \mathcal{K} \right\},
\end{equation}
where $\tau$ is an epigraphical variable.

\subsection{Model Discussion}\label{model_discussion}

In this section, we highlight some of the key characteristics of the models proposed above. First, both SO and RO have relatively complete recourse. That is, no matter what first-stage decisions we make, the second-stage problem always has a feasible solution. To verify this, consider a case where irrespective of the value of $z_{j}$'s, we set the value of $u_{j} = 0$ for all $j$ in the recourse function. This leads us to a feasible solution with no power generation and maximum load shed. 

Second, the proposed models can be further tightened based on some simple observations. To do so, we first compute the maximum flood level across all the scenarios for the flooded substations. Let us represent this value using parameter $W_i, \,\,\, \forall \, i \in \mathcal{I}_f$. Next, we observe that the model need not harden any substation to flood height that is higher than $W_i$. Therefore, in constraints (\ref{constraint_tighetening}), we can replace $H_i$ with min($H_i$, $W_i$), $\forall i \in \mathcal{I}_f$. Constraints \eqref{linking_1} and \eqref{linking_2} use the big-$M$ method. Here, we need to determine the smallest value for $M$ for each constraint. The smallest big-$M$ value for constraints \eqref{linking_1} and \eqref{linking_2} is given by $W_{\theta(j)}$ and min($H_{\theta(j)}$, $W_{\theta(j)}$) + 0.5, respectively. To verify this, recall the assumption that both the flood height and hardening level can only take non-negative integer values. Also, observe that 
$$-\min(H_{\theta(j)}, W_{\theta(j)}) \leq \Delta_{\theta(j)k} - x_{\theta(j)} \leq W_{\theta(j)}.$$ 
Now, in constraints \eqref{linking_1}, we need $\frac{\Delta_{\theta(j)k} - x_{\theta(j)}}{M} \leq 1$. The smallest value of $M$ that ensures this is $W_{\theta(j)}$. Similarly, in constraints \eqref{linking_2}, we need $\frac{1-2(\Delta_{\theta(j)k} - x_{\theta(j)})}{2M} \leq 1$. The smallest value of $M$ to achieve this is $\min(H_i, W_i) + 0.5$. 
Finally, we can also tighten constraints \eqref{phase_angle_constraint1} and \eqref{phase_angle_constraint2}. For both sets of constraints, the smallest value of $M$ is $F_r + 2 \pi B_r$.
 
Third, both SO and RO can be used for hardening decisions that will provide flood mitigation over a planning horizon. To see how, notice that the objective function in the SO computes the expected load shed for a single flood event. However, substation hardening in practice is done over the planning horizon that lasts multiple hurricane seasons and provides permanent protection for multiple flood events. Therefore, to help make hardening decisions that impact performance over multiple events, we first need to compute the disaster management costs due to load shedding over the multi-year planning horizon. To do so, let us assume that the expected number of hurricanes that the study region experiences during the planning horizon is $\gamma$. 
We further assume that during the planning horizon, the total recovery, economic and social costs are represented by the value of load loss of \$$\delta$/megawatt-hour. Finally, we assume that it takes $h$ hours to repair all the substations (and restore power to normal operation) starting immediately after a flood event. We believe this assumption is reasonable at the mitigation phase of the decision making process and avoids explicit and detailed modeling of the recovery process. Then, for a given investment budget, the substation hardening decisions that achieve the minimum expected total disaster management cost due to load shedding during the planning horizon is found by solving
\begin{equation}\label{total_cost}
    DM_{\mathcal{SO}} =  \gamma h \delta \, L_{\mathcal{S}\mathcal{O}}^*.
\end{equation}

In \eqref{total_cost}, 
the optimal substation hardening plan to minimize $DM_{\mathcal{SO}}$ is the same as the plan obtained by solving SO. This is because the $DM_{\mathcal{SO}}$ always equals the objective function of the SO multiplied by a positive constant. Therefore, irrespective of how the frequency of hurricanes, the restoration time, and the value of load loss change over time, the optimal substation hardening plan remains the same as the one obtained from SO. In this case, we assume that the probability distribution over the flood scenarios, and thus, the hurricanes causing them, does not significantly change over time. In practice, this is reasonable for planning horizons for which substation hardening is considered (5-10 years).
Similarly, the optimal substation hardening decision that minimizes the maximum total disaster management cost due to load shedding during the planning horizon is found by solving 
\begin{equation}\label{cost_bound}
DM_{\mathcal{RO}} = \gamma h \delta \, L_{\mathcal{R}\mathcal{O}}^*.
\end{equation} 
A key observation is that $DM_{\mathcal{RO}}$ provides an upper-bound for $DM_{\mathcal{SO}}$.
To understand why, note that the set of feasible first-stage solutions is same for both SO and RO. Further,  observe that both the models are bounded below with a minimum objective value of zero and have relatively complete recourse. Now, let $x_{\mathcal{R}}$ be a feasible solution to RO. Then,
\begin{subequations}
\allowdisplaybreaks
\begin{align}
    L_\mathcal{SO}^* & = \displaystyle{\minimize_{x\in \mathcal{X}} L_\mathcal{SO}(x)} \\     
    & \leq L_\mathcal{SO}(x_{\mathcal{R}}) \label{eq:robust_feasible}\\
    & = \sum_{k \in \mathcal{K}} p_k \, \mathcal{L}(x_{\mathcal{R}},k) \\
    & \leq \sum_{k \in \mathcal{K}} p_k \, \displaystyle{\maximize_{k \in \mathcal{K}}} \, \mathcal{L}(x_{\mathcal{R}},k) \label{eq:equal_load_shed} \\
    & = \displaystyle{\maximize_{k \in \mathcal{K}}} \,\, \mathcal{L}(x_{\mathcal{R}},k) \,\, \sum_{k \in \mathcal{K}} p_k\\
    & = \displaystyle{\maximize_{k \in \mathcal{K}}} \, \mathcal{L}(x_{\mathcal{R}},k) \\
    & = L_\mathcal{RO}(x_{\mathcal{R}}).
\end{align}
\end{subequations}
with Equation \eqref{eq:robust_feasible} holding at equality if and only if $x_{\mathcal{R}}$ is optimal for SO and Equation \eqref{eq:equal_load_shed} holding at equality if and only if load shed values are equal across all the scenarios. The above inequalities establish that the objective function value corresponding to any feasible solution to RO provides an upper bound on the optimal objective function value of SO. Since any optimal solution to RO is also feasible, it acts as a valid upper bound. In fact, it is the tightest upper bound that can be obtained in this manner. 

Finally, in the models proposed so far, we assume that we have a predetermined budget for substation hardening. One may however be interested in determining the optimal budget allocation for minimizing a risk measure over the disaster cost incurred due to both load shedding and substation hardening. The proposed models can easily be modified to find the optimal budget for substation hardening and corresponding hardening decisions. In the case of SO, this can be done by solving
\begin{equation}
TDM_{\mathcal{SO}} = \left\{ \displaystyle{\minimize  \, \, \omega \left (\sum_{k \in \mathcal{K}} p_k \, \mathcal{L}(x,k) \right) + \sum_{i \in \mathcal{I}_f} f_iy_i + v_ix_i}: x_{i} \leq H_i y_i, \,\, \forall  i \in \mathcal{I}_f \right\}, 
\end{equation}
where $\omega = \gamma h \delta$.
The value of $\sum_{i \in \mathcal{I}_f} f_iy_i + v_ix_i$ in the optimal solution represents the value of the optimal investment budget. Similarly, for RO, we compute the optimal investment budget by solving
\begin{equation}
TDM_{\mathcal{RO}} = \left\{ \displaystyle{\minimize  \, \, \omega \tau + \sum_{i \in \mathcal{I}_f} f_iy_i + v_ix_i}: \tau \geq \mathcal{L}(x,k) \,\, \forall \quad k \in \mathcal{K} \right\}.
\end{equation}

%% file: case.tex
\section{Case Study} \label{case_study}

In this section, using a case study for the Texas coastal region, we show how the proposed models can be used for power grid resilience decision making. The two main inputs to the proposed models are a set of scenarios that represent flood profiles for different hurricane types and the network parameters for the DC power flow model. To represent flood profiles, we use storm-surge maps developed by the National Oceanic and Atmospheric Administration (NOAA) \cite{MEOW}. For the electric grid, we use the ACTIVSg2000 dataset developed as part of an ARPA-E project \cite{activsg2000}. The details of each component are described in the following subsections. We further highlight that although we use the proposed models for storm surge-induced damages, they can be adopted for flooding of any kind as long as the corresponding flood scenarios are available. This could include scenarios for inland flooding as developed in \cite{kyoung_paper} and used for infrastructure resilience problems in \cite{SOUTO2022107545} and \cite{gizem_2022}. Lastly, to solve the various parameterizations of the proposed models discussed in this case study, we use the Gurobi solver with the barrier algorithm \cite{gurobi}. Within the solver, we set the MIP-gap threshold to 0.5 percent and limit the solve time to 6 hours. The model is solved on an Apple M1 pro machine with 16 GB of unified memory. 

\subsection{Flood Scenarios}

\begin{figure}
\captionsetup{width=0.45\textwidth}
\centering
\begin{minipage}{.5\textwidth}
    \centering
    \includegraphics[width=0.95\textwidth]{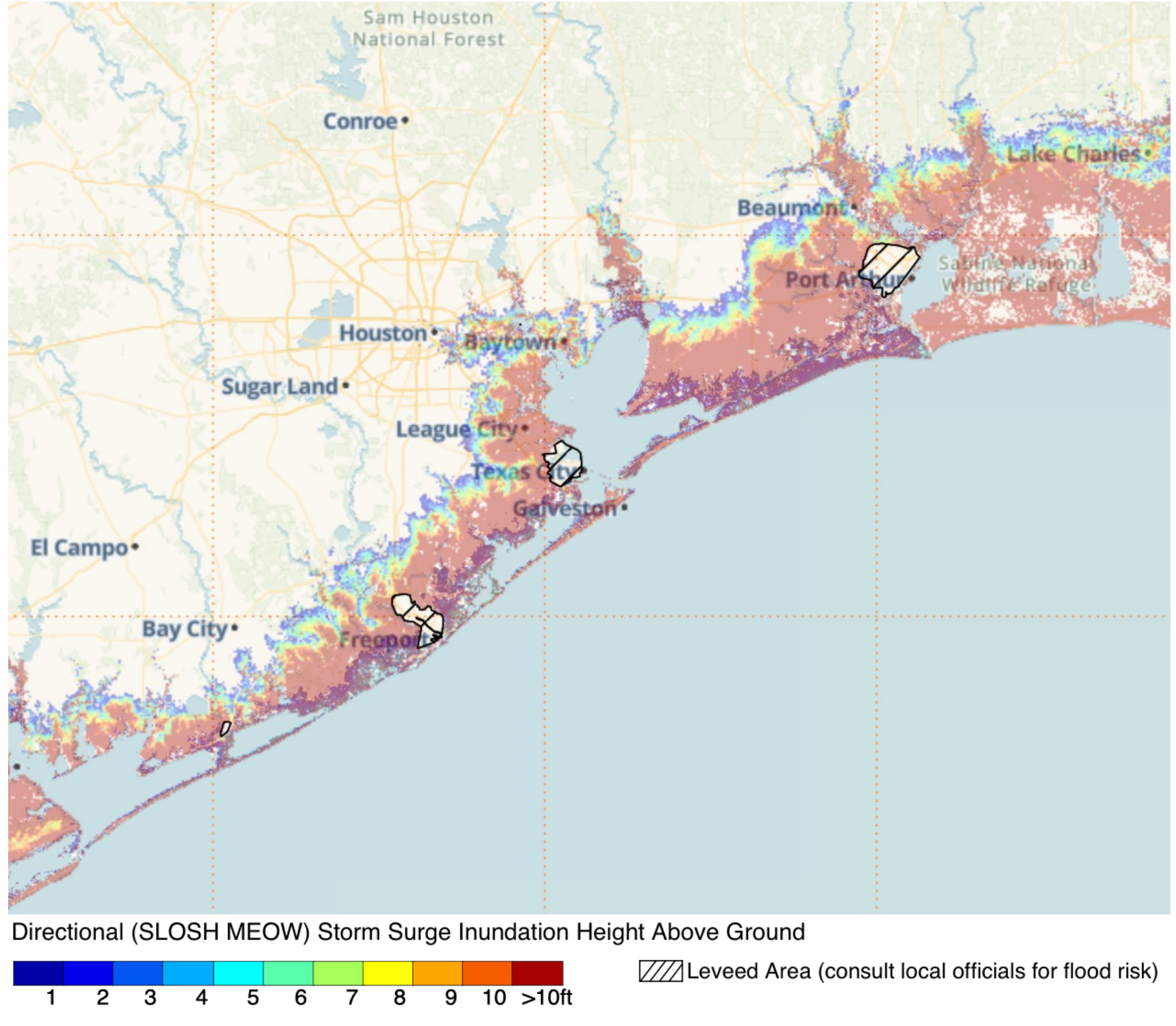}
    \captionof{figure}{A sample MEOW generated using category 5 storms approaching the Texas Gulf Coast in the north-west direction with a forward speed of 5 mph}
    \label{fig:MEOW}
\end{minipage}%
\begin{minipage}{.5\textwidth}
    \centering
    \captionof{table}{Grid characteristics before and after the electrically equivalent reduction was performed.}
    \begin{tabular}{|l|r|r|}
        \hline
        Grid Characteristic & Before & After \\
        \hline
        Substations (\#) & 1250 & 362 \\
        Buses (\#) & 2000 & 663 \\
        Transformers (\#) & 860 & 358 \\
        Transmission Lines (\#) & 2346 & 1151 \\
        Generators (\#) & 544 & 254 \\
        \hline
        Generation Capacity (GW) & 96.29 & 50.98 \\
        Load (GW) & 67.11 & 39.69 \\
        \hline
    \end{tabular}
    
    \label{tab:grid_before_and_after}
\end{minipage}
\end{figure}



We use the storm-surge maps developed by NOAA using the Sea, Lake, and Overland Surges from Hurricanes (SLOSH) model as flood scenarios. To generate these flood maps, SLOSH uses a simplified parametric wind field model that takes as input the following parameters: storm track, the radius of the maximum wind speeds, and the pressure differential between the storm's central pressure and the ambient pressure. The simulated wind fields are then used to compute surface stresses on the water beneath the hurricane. Finally, the induced stress on the surface of the water is used to determine the storm surge. For a detailed discussion on SLOSH, we refer to \cite{SLOSH}. Simulation studies developed using SLOSH have been extensively used to assist agencies like the Federal Emergency Management Administration (FEMA) and the U.S. Army Corps of Engineers (USACE).

In addition to real-time storm-surge guidance for imminent hurricanes, NOAA has developed two composite products, Maximum Envelopes of Water (MEOW) and Maximum of the MEOWs (MOM), to provide manageable datasets for medium to long-term hurricane evacuation planning. To develop these datasets, hurricane simulations with different combinations of intensity, forward speed, direction, and tide levels are run in parallel using SLOSH for the region of interest. 
Each run may yield different storm surge values for the same grid cell. A maximum overall such value is taken to represent the MEOW value of that grid cell. The same process is repeated for each grid cell within the study region to construct a MEOW map. The resolution of the grid cell is varied to balance accuracy with computation cost. It is finer in regions close to the coast and gets coarser as we move farther away in the ocean. MEOW maps are used to incorporate the uncertainties associated with a given forecast and eliminate the possibility that a critical storm track will be missed in which extreme storm surge values are generated. These maps are generated from several thousand SLOSH runs. In this study, we use the MEOW maps to represent flooding due to storm surge. An example MEOW map is shown in Figure \ref{fig:MEOW}. 

Within a MEOW map, it is possible that the water level for adjacent cells may come from different SLOSH runs of specific simulated storms.
Nevertheless, since these are the maximum water levels over multiple tracks, we can be assured that if we have hardened a substation for a particular MEOW map, it will provide resilience towards flooding for any of the parallel runs that constitute the MEOW map. Arguably, MEOW maps are still better at representing flood uncertainty than other scenario generation methods where flooding at different nodes within a network is considered independent of each other. Moreover, these MEOW maps have been considered as scenarios in different stochastic optimization models for patient evacuation \cite{kyoung_paper} and grid hardening \cite{shukla_scenario-based_2021}. 

The original MEOW map dataset for the Texas coastal region comprises 192 flood maps which are constructed using eight different storm directions (west-south-west, west, west-north-west, north-west, north-north-west, north, north-north-east, and north-east), six different intensity categories (0-5), and four different, forward speeds (5,10,15, and 25 mph). To demonstrate the usefulness of the proposed approach with a computationally tractable use case, we reduce the size of the problem by eliminating a subset of less severe scenarios. We first drop all the MEOW maps corresponding to four directions (west-south-west, north, north-north-east, and northeast) as hurricanes belonging to these categories do not cause significant flooding in the Texas Gulf Coast. We also drop the MEOW maps corresponding to category 0-4 hurricanes. The storms belonging to category 5 are more intense versions of these storms. Hence, the model ends up recommending decisions to prepare for worst-case situations (and thus implicitly prepares for category 0-4 hurricanes as well). Moreover, as discussed before, we  assume that in the mitigation phase, the decision-maker does not have information on any specific storm. Therefore, to model the uncertainty in the mitigation phase, we assume that all the remaining MEOW maps are representative of the flooding scenarios. In our case, they are considered equally likely for the stochastic model. This is based on the premise that the larger set of simulations that produced the MEOW maps were sampled according to some underlying distribution and therefore already reflects the underlying characteristics of the distribution implicitly used by NOAA in their development of the MEOW product. Furthermore, these MEOW maps provide us with the storm-surge flood height above ground at each of the substations (and thus for buses within) in the power grid network. In the proposed models, this is represented by $\Delta_{ik}$; the level of flooding at substation $i$ in scenario $k$.

\subsection{Power Grid}

\begin{figure}
    \centering
    \includegraphics[width=0.9\textwidth]{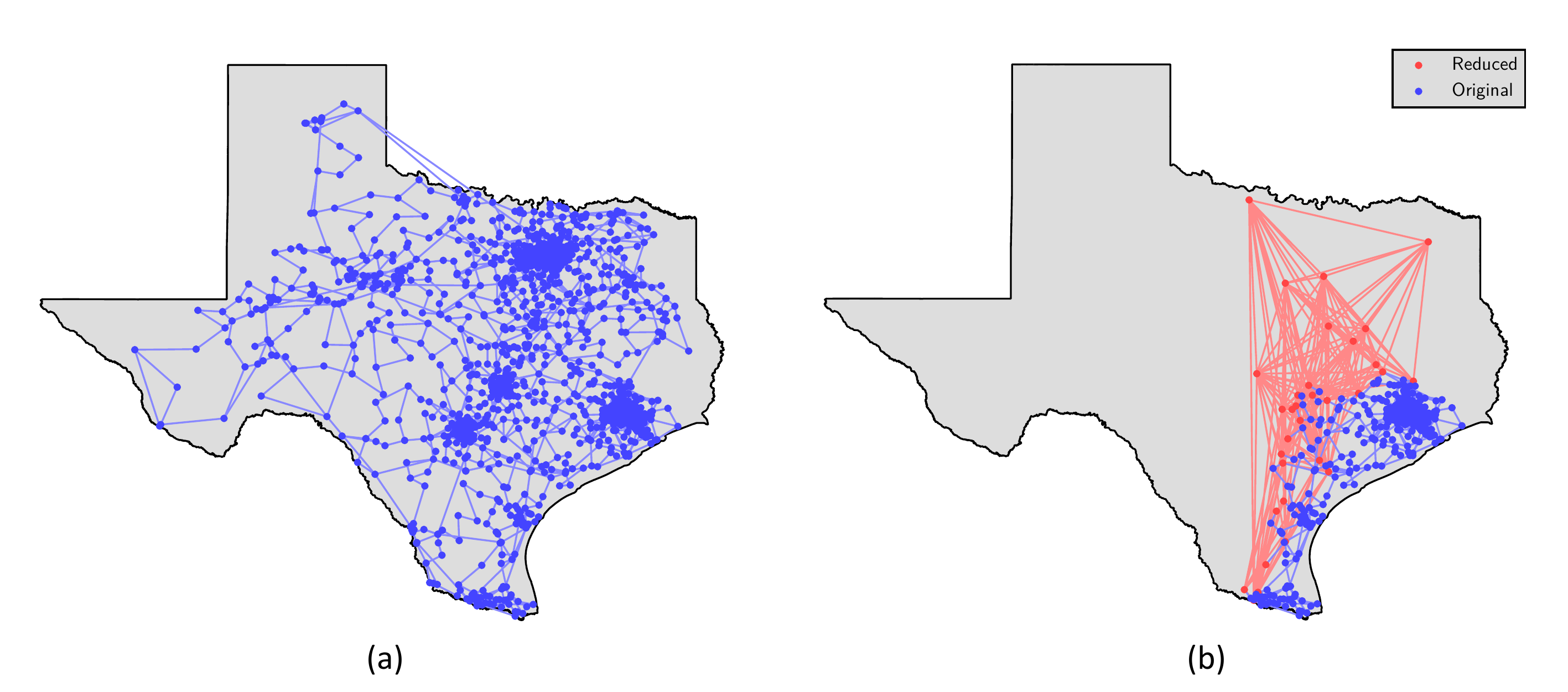}
    \caption{The figure shows (a) ACTIVSg2000 Synthetic Grid for Texas, and (b) the reduced grid obtained after performing the network reduction. The red elements represent the new nodes and branches that were introduced as the artifacts of the reduction procedure to maintain equivalence in the grid characteristics}
    \label{fig:grid_reduction}
\end{figure}

To model the power grid for the state of Texas, we use the synthetic grid called ACTIVSg2000 which contains 2000 buses (within 1250 substations) and 3206 branches. The grid, though synthetic, is designed such that it maintains statistical similarities with the actual Texas grid. We make two further modifications to the grid instance to make it computationally tractable while also considering the coastal part of the grid which is affected due to to storm-surge and thus is the focus of the study. First, we perform a network reduction on the original grid instance using the electrical equivalent (EEQV) feature in {PSS\textregistered E} to focus on the grid components subject to storm surge flooding. The reduction is such that the buses in the inland region that are not exposed to storm-surge induced flooding are aggregated within a much smaller set of nodes. The part of the grid that is in close proximity to the Texas Gulf Coast, and therefore is prone to flooding due to storm-surge, is retained almost as is. The effect of the network reduction is detailed in Table~\ref{tab:grid_before_and_after}. The topological changes are visualized in Figure \ref{fig:grid_reduction}. Second, we alter the locations of some of the substations. This is because, in the original dataset, some of the substations are placed in the middle of a water body and are thus flooded by default. To address this, we remap the coordinates of the 1250 substations in the dataset with the coordinates of substations obtained from the Homeland Infrastructure Foundation-Level Data (HIFLD) Electric Substations dataset \cite{HIFLD_Electric_Substations}. The HFILD dataset contains information about real-world substations across the U.S. The remapping is done by solving an optimization problem that minimizes the total displacement due to remapping. Note that this process does not change the power grid's electrical structure and makes it more realistic using the real-world substation locations (closer to the actual Texas grid) and capturing their real-world flood risks via MEOW-based flood scenarios. Lastly, the fixed cost and the variable cost for substation hardening are assumed to be \$25,000 and \$100,000 per foot, respectively. These values are derived from various utility reports.  


%% file: results.tex
\section{Results and Discussion} \label{results}

In this section, we first determine the expected value of perfect information from a model that can produce perfect forecasts. In the same subsection, we compute the value of the stochastic solution for the different budget levels. Next, in subsection \ref{opti_budget} we show how the proposed model can be used to determine the optimal investment budget for substation hardening. 

Due to dynamically changing climate conditions and ocean temperature, the probabilities of different types of hurricanes can change over time. In subsection \ref{uncertain prob}, we show when a decision-maker can take advantage of using solutions from RO to hedge against this uncertainty by paying a relatively insignificant premium. Finally, the last subsection is dedicated to the analysis of the distribution of load shed across scenarios for the solutions obtained from SO and RO.

\subsection{The Values of Stochastic Solution and Perfect Information}
\label{vssevpi}

In the proposed two-stage decision-making models, both the number of variables and constraints grow linearly with the number of scenarios; thus making it computationally challenging. In that case, instead of solving SO for large grid instances with many scenarios, one may be interested in solving simpler versions of the problem. One approach could be to reduce the size of the problem by constructing a single scenario problem where the flood height at each substation is the average of the flood height across all the scenarios. Another way could be to solve the problem for each scenario individually. The first-stage solutions thus obtained can then be analyzed and potentially combined using some heuristic rule. In this section, we analyze the quality of the solutions we  get using such approaches. To do so, we use two widely known concepts in the stochastic programming literature: the value of the stochastic solution and the expected value of perfect information. 

\begin{figure}
\captionsetup{width=0.45\textwidth}
\centering
\begin{minipage}{.5\textwidth}
    \centering
    \includegraphics[scale=0.7]{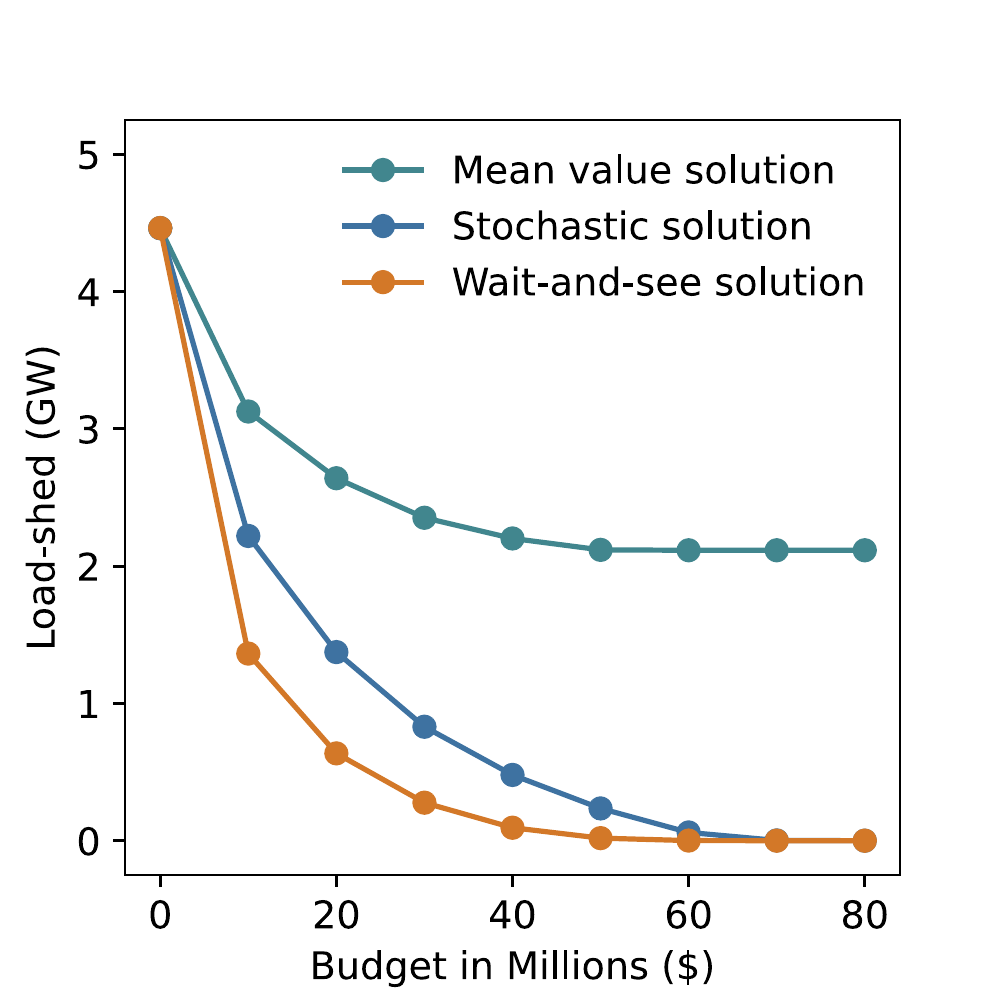}
    \captionof{figure}{The graph of the expected load shed values for the expected value solution (green), the stochastic solution (blue), and the wait-and-see solution (orange), as a function of the budget for substation hardening}
    \label{fig:16_scenarios_bounds}
\end{minipage}%
\begin{minipage}{.45\textwidth}
    \centering
    \includegraphics[scale=0.7]{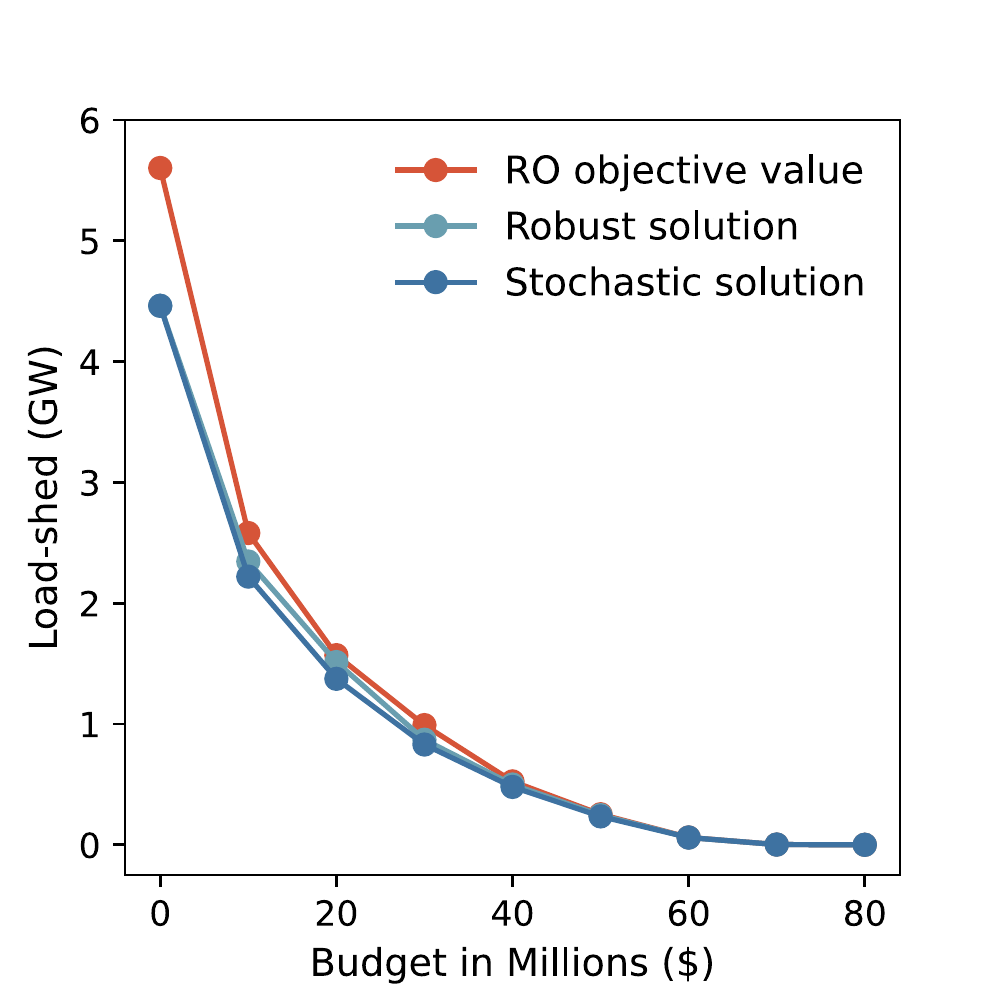}
    \captionof{figure}{The objective function value for RO (i.e., the optimal maximum load shed) and the expected load shed values for the robust and the stochastic solutions as a function of budget for substation hardening}
    \label{fig:robust_bound}
\end{minipage}
\end{figure}

In Figure \ref{fig:16_scenarios_bounds}, we plot the value of $L_{\mathcal{SO}}^*$ as a function of the investment budget. To determine the budget levels on which the parametric study should be performed, we first compute the minimum budget such that $L_{\mathcal{SO}}^* = 0$. This is computed by solving a slightly modified version of SO. Specifically, we first remove constraint \eqref{budget_constraint} from the formulation and replace the objective function with the minimization of the substation hardening expenditure (i.e., the left-hand-side of \eqref{budget_constraint}). We also force full satisfaction of demand by replacing constraint \eqref{supply_demand} with $s_{j} = D_j, \,\,\, \forall \, j \in \mathcal{J}.$ The optimal value to this modified version of SO is in turn the minimum hardening budget required for zero load shed. For the parameters assumed in this case study, the corresponding minimum budget turns out to be \$71.35M. Any additional budget beyond this will not improve the objective function value (load shed) and therefore the corresponding optimal solution. Using \$71.35M as the reference, in Figure \ref{fig:16_scenarios_bounds}, we increase the budget from \$0M, in increments of \$10M, until the value exceeds \$71.35M. (We note that the same budget values can be used for the RO parametric study. This is because the minimum budget required to achieve the expected load shed of zero is the same as what is required to achieve zero load shed in all the scenarios.)



In addition to computing $L_{\mathcal{SO}}^*$ for different budget levels, we also compute lower and upper bounds on this value as shown in Figure \ref{fig:16_scenarios_bounds}. To compute the upper bound, let us first consider a single-scenario problem called the expected value problem which is defined as
\begin{equation}\label{ev}
L_{\mathcal{EV}}^* = \displaystyle{\minimize_{x\in \mathcal{X}} \, \mathcal{L}(x,\bar{k})}.    
\end{equation}
Here, $\bar{k}$ represents a scenario where the flood level at each substation is the mean of the flood heights at that particular substation across all the scenarios. The optimal solution to this problem is called the expected value solution or the mean value solution, henceforth represented by $\bar{x}$. We note that the problem in \eqref{ev} is a single scenario problem and therefore much smaller in size. However, the reduction of scenarios leads to the loss of information about the substations that flood together. Next, to evaluate the quality of the first-stage substation hardening decisions obtained by solving \eqref{ev}, we compute the expected load shed across all the scenarios by fixing the first-stage decisions to $\bar{x}$, denoted by $L_{\mathcal{SO}}(\bar{x})$. This value serves as an upper bound on $L_{\mathcal{SO}}^*$. For a detailed explanation on this, we refer to \cite{BirgLouv97}. We also observe in Figure \ref{fig:16_scenarios_bounds} that the difference between $L_{\mathcal{SO}}(\bar{x})$ and $L_{\mathcal{SO}}^*$ increases with an increase in the budget. We refer to this difference as the value of the stochastic solution (VSS) as it represents the value of using a scenario-based representation of the uncertainty as opposed to the average flood values, all calculated within the SO framework. We further highlight that when $\bar{x}$ is implemented, the load shed does not strictly decrease with an increase in the budget beyond \$50M. This is expected because, for the mean scenario, the model obtains zero load shed with \$50.15M. As discussed before, we know that it takes a minimum of \$71.35M to achieve zero load shed across all the scenarios. However, once the model achieves zero load shed in the mean scenario, there is no incentive to use additional resources. This also shows that the value of using SO over the expected value problem increases with increases in the investment budget.

To compute a lower bound on $L_{\mathcal{SO}}^*$ for a given budget level, we assume that the decision-maker has access to perfect information about the flood levels, and therefore can better prepare for each scenario (i.e., in a way, fine-tuning the mitigation plan according to each scenario). That is, perfect information allows the decision maker to make possibly different substation hardening decisions in each scenario to minimize the load shed in that particular scenario. These solutions are referred to as wait-and-see solutions. In this case, we compute
\begin{equation}
    L_{\mathcal{WS}}^* = \sum_{k \in \mathcal{K}} p_k \, \left( 
\displaystyle{\minimize_{x\in \mathcal{X}}\mathcal{L}(x,k)}
    \right),
\end{equation}
where the first-stage decisions are scenario-dependent and the value $L_{\mathcal{WS}}^*$ is referred to as the wait-and-see bound. Due to the scenario-specific mitigation decisions, $L_{\mathcal{WS}}^*$ provides a lower bound on $L_{\mathcal{SO}}^*$. The difference in the values $L_{\mathcal{SO}}^*$ and $L_{\mathcal{WS}}^*$ is referred to as the expected value of perfect information. It represents the maximum value a decision-maker would be willing to pay in exchange for complete and accurate information about the uncertainty. 

The key point that we want to emphasize is that unless the flood model can make perfect predictions, which is usually not the case with the weather models, then not accounting for uncertainty and using just point estimates or mean values of the flood forecasts can lead to significantly inferior decisions. This is evident from the VSS. In fact, as it turns out in this case, the first-stage decisions obtained from SO, even with not-so-perfect forecasts, lead to a load shed performance that is close to what we would get from using a flood model that offers perfect prediction. To put it another way, accounting for flood uncertainty, even with a small number of scenarios, can help reduce the burden of getting perfect information on the decision-maker without making a significant compromise on the performance.

We also note that all bounds converge to the same expected value at the zero budget. This is because, no matter how well we represent the uncertainty or how good the predictions are, if we do not have any resources to use towards mitigation in the first stage, we cannot prevent load shed in the second stage with no protection towards flooding. Then, the expected value of the load shed is  only a function of the second-stage decisions (i.e., the best power flow the grid can deliver with flooded substations). Moreover, at a sufficiently high budget value, both $L_{\mathcal{SO}}^*$ and $L_{\mathcal{WS}}^*$ converge to zero. This is expected because, despite the fact that we have poor predictions or poor uncertainty representation, we can still prevent any load shed in all the scenarios if we have enough resources to harden all substations to any desirable extent.

Performing analysis with the budget level as a parameter, as described in this subsection, requires repeatedly solving SO with different parameters and can be time-consuming. To address this, the property that SO has relatively complete recourse is exploited to warm start the optimization solver and improve the solution time. As was stated in subsection \ref{model_discussion}, we can heuristically generate an initial feasible solution with a hundred percent load shedding for an investment budget value of zero. Once we get an optimal solution corresponding to budget level zero, we use it to warm start the solver for the next higher budget level. The process is repeated to generate high-quality feasible solutions for the next budget level. We further note that the aforementioned approach for warm-starting the solver is also applicable in the case of RO.

\subsection{Determining optimal budget for substation hardening}\label{opti_budget}
Subsection \ref{vssevpi} focuses on demonstrating how SO is used for resilience decision making for a given investment budget for substation hardening. The models can alternatively be used to decide the optimal value of the investment budget that minimizes the expected total disaster management cost over a multi-year planning horizon. To demonstrate this, we solve $TDM_{\mathcal{SO}}$ as described in Section \ref{model_discussion} assuming that, on average, 10 hurricanes hit the Texas Gulf Coast during the planning horizon. Furthermore, for sensitivity analysis, we consider three different values of restoration time: 6, 24, and 48 hours. Similarly, for the value of load loss, we consider 5 different values: \$250, \$500, \$1000, \$3000, and \$5000 per MWh. 
\begin{figure}[ht]
    \centering
    \includegraphics[scale=0.8]{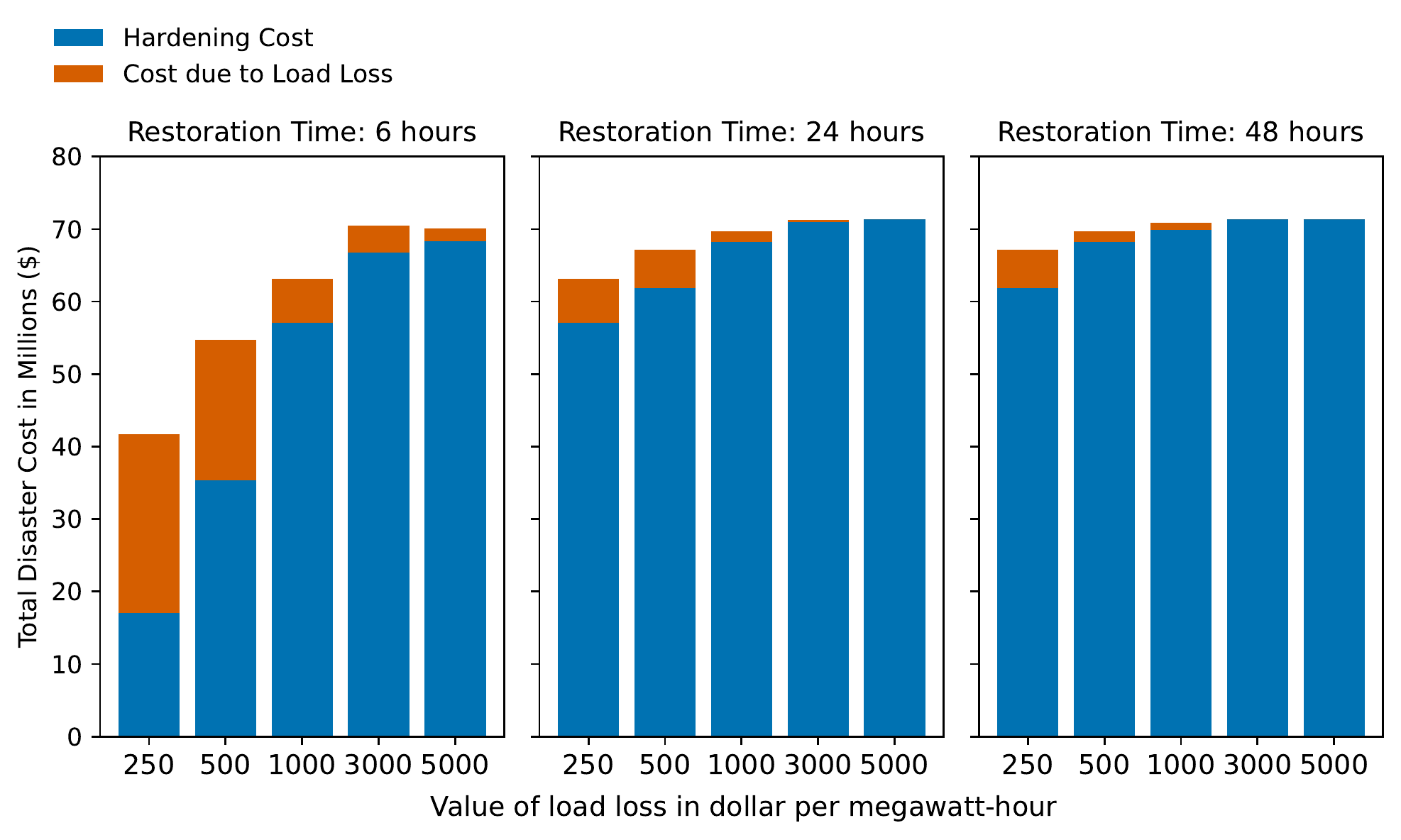}
    \caption{Total disaster management cost as a function of the value of load loss for different restoration times}
    \label{fig:voll_plot}
\end{figure}

The values of the expected total disaster management cost for  different combinations of restoration times and the value of load loss are plotted in Figure \ref{fig:voll_plot}. As we see in the figure, when the value of load loss is low and the restoration time is short (\$250 and 6 hours, respectively), the total disaster management cost is relatively low. Furthermore, the model recommends investing only a quarter of the budget required to achieve zero load shed for substation hardening. This is because the cost associated with losing power is quite low and the outage is restored relatively quickly. On the other hand, when both the load loss value and restoration time are high (\$5000 and 48 hours, respectively), the model recommends investing \$71.35M (equal to the investment required to achieve zero load shed) for substation hardening, avoiding any costs due to load loss. This is because, for the chosen values, the costs associated with power loss are quite high and it is better to make significant investments (in fact, everything possible) in substation hardening such that there is no loss of power. For all other combinations in between, both the total expected disaster management cost and its composition vary. Based on a study undertaken by the Electric Reliability Council of Texas, the value of load loss was determined to be around \$6000 per MWh for Texas \cite{ercot_voll}. If we take that as given, the results in Figure \ref{fig:voll_plot} suggest that we must make investments to achieve close to zero load shed even if restoration time is as short as 6 hours.

\begin{figure}[h]
    \centering
    \includegraphics[scale=0.7]{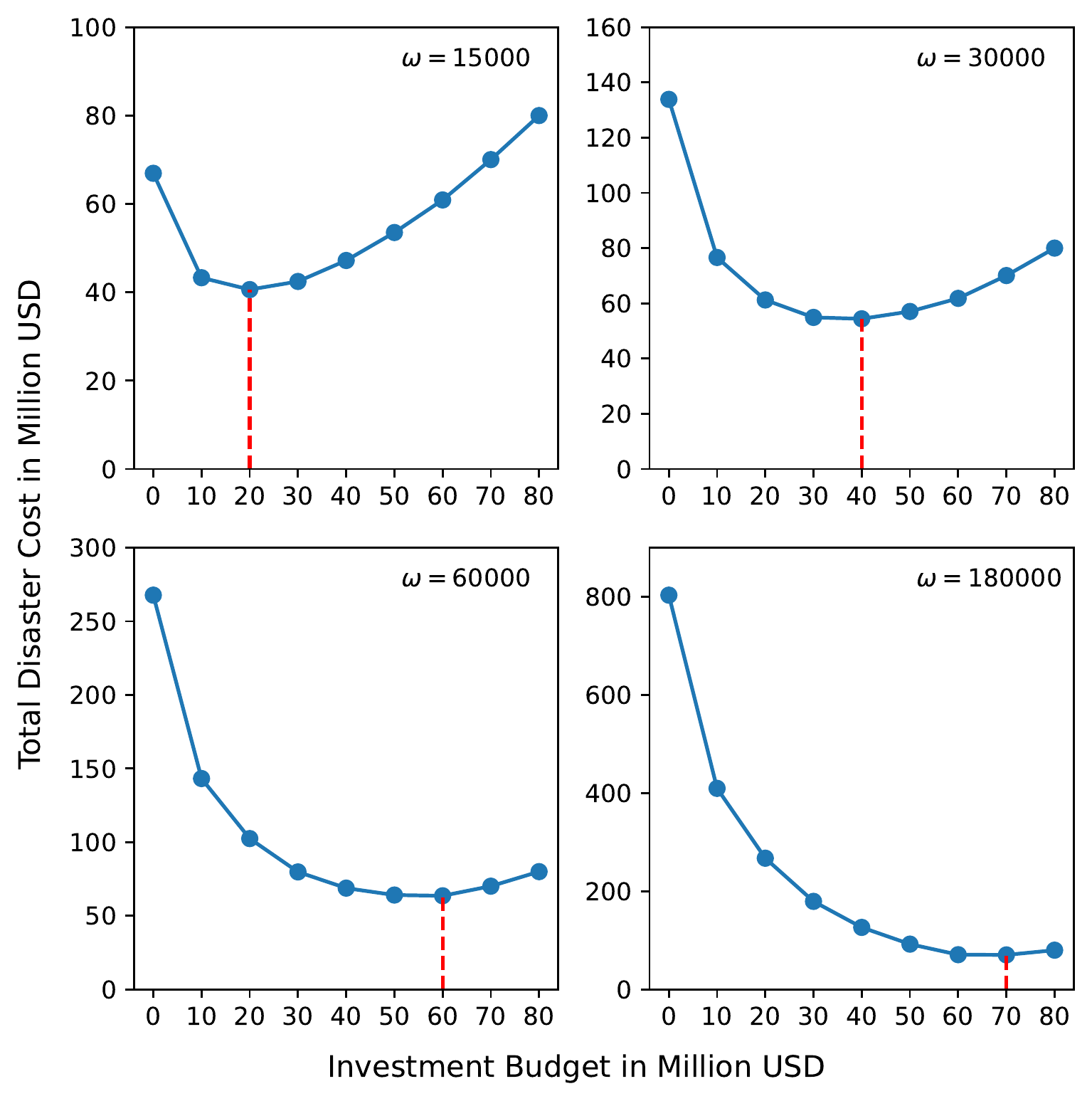}
    \caption{Total disaster management cost as a function of budget for different values of $\omega$}
    \label{fig:voll_quick}
\end{figure}

We further observe in Figure \ref{fig:voll_plot} that the solution corresponding to the value of load loss of \$1000 per MWh and restoration time of 6 hours is the same as the solution with the value of load loss of \$250 per MWh and a restoration time of 24 hours. This is expected because in $TDM_{\mathcal{SO}}$, both sets of parameters lead to the same optimization problem ($\omega$ is the same). It is further apparent that the optimal investment budget for substation hardening increases monotonically from 0 to \$71.35M with the increase in the value of $\omega$. Therefore, for any investment budget value between 0 and \$71.35M, there exists a unique $\omega$ for which the corresponding budget is optimal. We can use this insight to quickly approximate the optimal investment budget for any combination of the value of load loss, restoration time, and the average number of hurricanes that may hit the region of study during the planning horizon. To do so, we use only the values of $L_{\mathcal{SO}}^*$ for $I \in \{0,10,20...80\}$ as computed in Section \ref{vssevpi}. For any given value of $\gamma, h,$ and $\delta$ for which the optimal investment budget needs to be approximated, we compute the value of $DM_{\mathcal{SO}} + I$ for each $I \in \{0,10,20...80\}$. The value of $I$ for which $DM_{\mathcal{SO}} + I$ is the smallest is the best approximation for the optimal investment budget for a chosen value $\omega$ (calculated from $\gamma, h,$ and $\delta$). In Figure \ref{fig:voll_quick}, we show the value of $DM_{\mathcal{SO}} + I$ for $I \in \{0,10,20...80\}$ for different values of $\omega$. The depicted values of $\omega$ are reasonable in the sense that they can be derived from the $\gamma, h,$ and $\delta$ used in Figure \ref{fig:voll_plot}. For example, we notice in Figure \ref{fig:voll_plot} that when $\gamma = 10, h= 6,$ and  $\delta = 250$ ($\omega = 15000$), the optimal investment budget is \$17.05M. An approximate of this value can be quickly inferred by looking at Figure \ref{fig:voll_quick} for the value of $\omega = 15000$. For that value, the top-left curve achieves its minimum at \$20M which is closest to \$17.05M in the set $\{0,10,20...80\}$.

\subsection{Optimization in the face of uncertain probabilities}\label{uncertain prob}

While framing the power-grid resilience decision-making problem as a two-stage stochastic program, we assume that the probability distribution over the scenarios is known. However, the probabilities that we assign to each scenario need not be constant with time. Changing climate oscillations and ocean temperatures routinely affect these probabilities. This is reflected in NOAA's annual hurricane season prediction categories: normal, above normal, and below normal. In this case, if probabilities over scenarios change as compared to what we planned for, the expected performance may deteriorate. To hedge against this, we recommend solving RO and comparing the expected performance of the corresponding decisions. If the difference in the expected performance of decisions recommended by RO and SO is not significant, we advise adopting decisions as recommended by RO. In this way, irrespective of how the probabilities evolve with time, we know that the expected total cost due to the load loss will be less than or equal to the bound obtained in Section \ref{model_discussion} if the optimal first-stage decisions as recommended by RO are adopted. This can be confirmed for the parameters assumed in the case study through Figure \ref{fig:robust_bound}. The robust solution refers to the value of the expected load shed when the first-stage hardening decisions as recommended by RO are adopted. Since these decisions are not necessarily optimal for the assumed probability distribution, they lead to a higher expected load shed as compared to the stochastic solution. However, for the RO first-stage decisions, the maximum load shed in any scenario is capped by $\tau$ as represented by the red curve. In this case, we observe that the difference in the expected value of performance is almost trivial for investment budget values of \$40M and above. Therefore, in those cases, it makes sense to adopt decisions recommended by RO as opposed to what we get from SO to hedge against the change in probabilities due to factors like ocean temperature and climate oscillations. In cases when the difference between the expected performance of SO and RO is significant, the decisions depend on the risk preference of the decision-maker. 

\subsection{SO vs RO: Analysis of the load shed distribution}

\begin{figure}[h]
    \centering
    \includegraphics[scale=0.4]{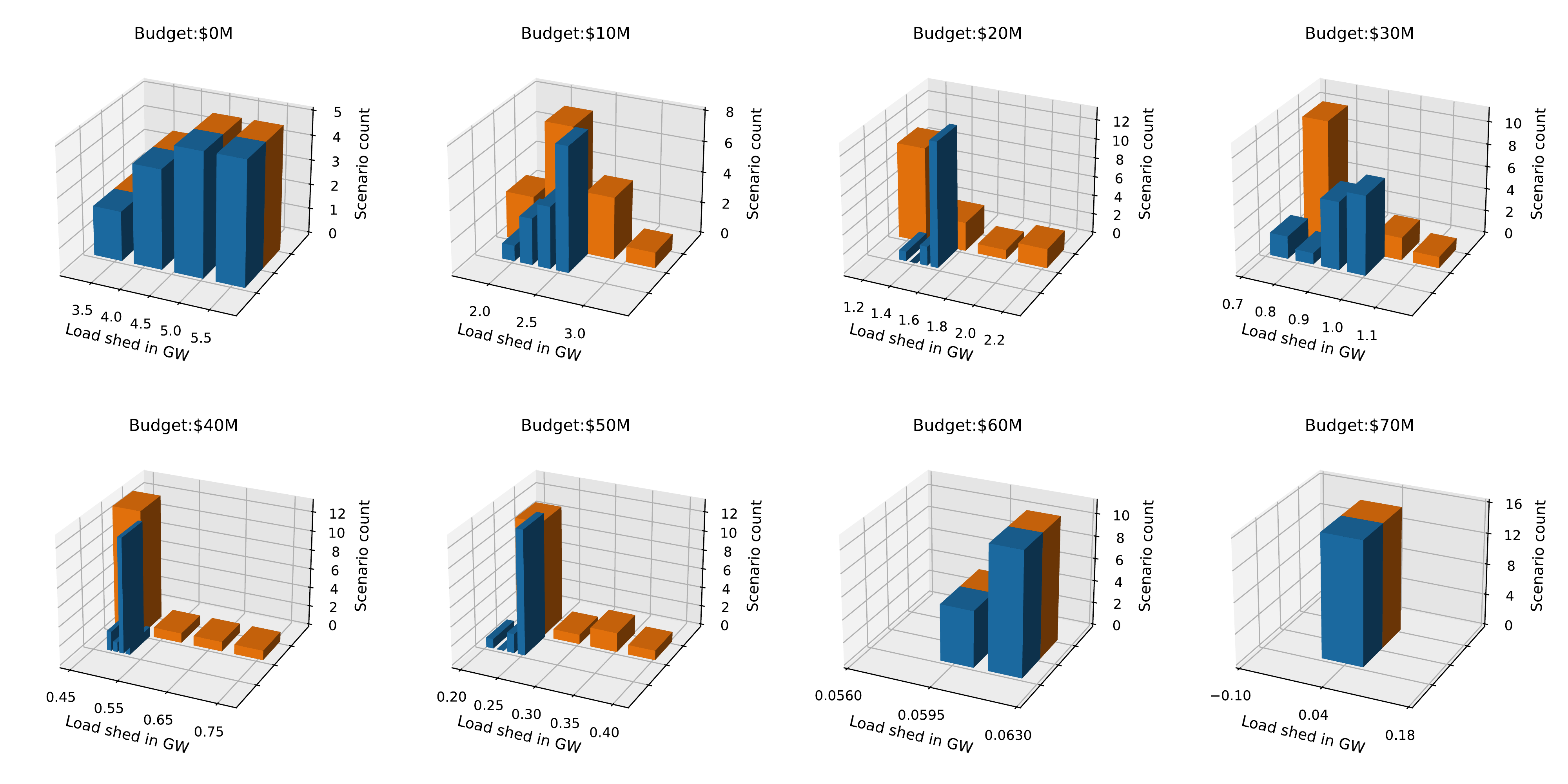}
    \caption{The histograms of load shed across scenarios for the SO and RO solutions at different budget levels. Notice that the x-axis for each sub-plot has a different scale for better depiction.}
    \label{fig:mixed_distribution}
\end{figure}

We conclude the discussion with an analysis of the distribution of the load shed across scenarios for both SO and RO. The load shed in each scenario for both models is represented by the value of the recourse function corresponding to the optimal solution. These values are used to construct the corresponding histogram for both SO and RO at different budget levels as shown in Figure \ref{fig:mixed_distribution}. 
As expected, the histograms shift to the left with the increase in the budget for both SO and RO. We observe that the histograms for both SO and RO coincide when the investment budget is \$0M. This is expected because there is no hardening done in either model. Consequently, the load shed in each scenario is identical. Moreover, using Figure \ref{fig:robust_bound} and Figure \ref{fig:mixed_distribution}, we observe that the robust solutions provide an inferior performance in expectation but the RO load shed remains relatively stable across scenarios as compared to SO. We  also conclude that for investment budget values of more than \$40M, the robust solutions offer much better performance against extreme scenarios while also offering good expected value performance. Therefore, in this case, it is reasonable to implement RO decisions for budget values beyond \$40M. In this way, Figures \ref{fig:robust_bound} and \ref{fig:mixed_distribution} can be used together to understand the behavior of both SO and RO in expectation and across all the individual scenarios. 

\section{Conclusions}\label{conclusions}

In this study, we propose an integrated framework supported by two scenario-based optimization models for power grid resilience decision making against extreme flooding events. The models recommend an optimal substation hardening plan by integrating the predictions generated from a state-of-the-art hydrological model with a DC optimal power flow model. While doing so, we account for uncertainty in hurricane predictions using a scenario-based representation. Furthermore, using a case study for the state of Texas, specifically the coastal region prone to storm surge flooding, we demonstrate how the proposed models can together be used to address a wide variety of insight-seeking questions related to power grid resilience decision making. Specifically, we show that using a scenario-based representation of flood uncertainty can offer significant value over mean flood forecasts. We also explain the advantages of using flood maps generated from physics-based models as opposed to other scenario generation methods popular in the literature. Furthermore, we show how can we estimate the expected value of perfect information from near-perfect flood forecasts. For the case study developed in the paper, we observe that by using a scenario-based representation of uncertainty, the decision-makers can reduce their burden of having access to perfect forecasts. In addition to quantifying the value of using flood scenarios, we further show how we can use the proposed two-stage framework to determine the optimal investment budget for substation hardening. Lastly, we explain how we can use the two-stage robust optimization model for power-grid resilience decision making when information about the probability distribution over the flood scenarios is unavailable.

For future research, we suggest four directions. First is the development of scenarios that can consider precipitation-induced inland flooding in addition to storm-surge. Second is the development of methods to account for equity while making substation-hardening decisions. Third is developing models that take into account preparedness measures while making longer-term mitigation decisions, leading to three-stage optimization models. Fourth is the development of decomposition techniques to solve such models in a reasonable time. These challenges form the basis of our ongoing research.